\documentclass[12 pt, english]{article}
\usepackage[T1]{fontenc}
\usepackage{ulem}
\usepackage{url}
\usepackage[utf8]{luainputenc}
\usepackage{geometry}
\usepackage[pdftex]{graphicx}
\usepackage{amstext}
\usepackage{amsthm}
\usepackage{amssymb}
\usepackage{amsmath}
\usepackage[T1,T2A]{fontenc}
\usepackage[utf8]{inputenc}
\usepackage[main=english]{babel}
\usepackage{setspace}
\usepackage{esint}
\usepackage{comment}
\usepackage{babel}
\usepackage{float}
\usepackage{amsfonts}
\usepackage{fullpage} 
\usepackage{parskip} 
\usepackage{tikz} 
\usepackage{indentfirst}
\renewcommand{\:}{\colon}

\newcommand{\floor}[1]{\left\lfloor #1 \right\rfloor}

\def\Area{\operatorname{Area}}
\def\diam{\operatorname{diam}}
\def\dim{\operatorname{dim}}

\def\spt{\operatorname{spt}}
\title{Besicovitch--type inequality for closed geodesics on 2-dimensional spheres}
\author{Talant Talipov}
\date{}
\begin{document}
\thispagestyle{empty}
\theoremstyle{definition}
\newtheorem{theorem}{Theorem}[section]
\newtheorem{lemma}{Lemma}[section]
\newtheorem{remark}{Remark}[section]
\newtheorem{corollary}{Corollary}[section]
\newtheorem{example}{Example}[section]
\newtheorem{definition}{Definition}[section]
\newtheorem{proposition}{Proposition}[section]
\newtheorem{conjecture}{Conjecture}[section]
\sloppy

\maketitle

\begin{abstract}
    We prove the existence of a constant $C > 0$ such that for any Riemannian metric $g$ on a 2-dimensional sphere $S^2$, there exist two distinct closed geodesics with lengths $L_{1}$ and $L_{2}$ satisfying $L_{1} L_{2} \leq C \cdot \Area(S^2, g)$.
\end{abstract}

\section{Introduction}
\hspace*{1cm}In this paper, we address a question posed by Y. Liokumovich in~\cite{LectureNotes}, which asks whether there exists a constant $C > 0$ such that for any 2-dimensional Riemannian sphere $(S^2, g)$, there exist two distinct closed geodesics with lengths $L_1$ and $L_2$ satisfying
\begin{equation*}
    L_1 L_2 \leq C \cdot \Area(S^2, g).
\end{equation*}
We prove that this conjecture is true.\\
\hspace*{1cm}The study of bounds on the lengths of closed geodesics has a rich history. In~\cite{Gromov}, M. Gromov established the systolic inequality, which bounds the length of the shortest non-contractible loop on a essential Riemannian manifold relative to its volume. Optimal constants for the systolic inequality have been determined for specific surfaces: the 2-dimensional torus by L. C. Lowner, the real projective plane by P. M. Pu (see~\cite{Pu}), and the Klein bottle by C. B. Croke and M. Katz (see~\cite{CrokeKatz}). However, the sphere is not an essential manifold.  The first bound on the length of the shortest closed geodesic on a Riemannian 2-sphere in terms of area was established by C. B. Croke~\cite{Croke}. This result was later improved, with the most recent enhancement made by R. Rotman, who provided the best currently known constant in~\cite{Rotman}. Furthermore, in \cite{LioNabRot}, Y. Liokumovich, A. Nabutovsky, and R. Rotman proved the existence of three distinct simple closed geodesics on a 2-dimensional Riemannian sphere, with their lengths bounded by the diameter of the sphere. Notably, it is known that the length of the second shortest closed geodesic cannot be bounded purely in terms of the square root of the area, as demonstrated by the example of long ellipsoids. Our result can be seen as a spherical analog of the Besicovitch inequality, which provides a bound on the product of distances between opposite sides of a Riemannian square in terms of its area.\\
\hspace*{1cm}We state two conjectured Besicovitch-type inequalities for the higher-dimensional case. Let $n > 2$.
\begin{conjecture}
There exists a constant $C(n) > 0$, depending on the dimension $n$, such that for any $n$-dimensional Riemannian sphere $(S^n, g)$, there exist $n$ distinct closed geodesics with lengths $L_1, \ldots, L_n$ satisfying
\begin{equation*}
    L_1 \cdot L_2 \cdot \ldots \cdot L_n \leq C(n) \cdot \operatorname{Vol}(S^n, g).
\end{equation*}
\end{conjecture}
\begin{conjecture}
There exists a constant $C(n) > 0$, depending on the dimension $n$, such that for any $n$-dimensional Riemannian sphere $(S^n, g)$, there exist a closed geodesic $\gamma$ and a smooth closed embedded minimal hypersurface $N^{n-1}$ satisfying
\begin{equation*}
    \operatorname{Length}(\gamma) \cdot \operatorname{Vol}_{n-1}(N^{n-1}) \leq C(n) \cdot \operatorname{Vol}(S^n, g).
\end{equation*}
\end{conjecture}
\hspace*{1cm}\textbf{Outline.} The paper is structured as follows. In section~\ref{sec-2}, we present a general proof of the statement using the Almgren--Pitts min-max theory. In section~\ref{sec-3}, we treat two specific cases based on the type of the shortest closed geodesic: one case corresponds to a simple closed geodesic with Morse index 1, and the other to a figure‐eight geodesic. In these cases, the multiplicative constant in our bound is lower than in the general case, and additional geometric properties of the second geodesic are established.\\
\hspace*{1cm}In the first case, referred to as \textit{long sphere}, we assume that both two connected components of the complement of the shortest closed geodesic have diameters sufficiently larger $\frac{\operatorname{Area}(S^2)}{L_1}$, where $L_1$ is the length of the shortest closed geodesic. By applying the Birkhoff curve shortening process to a suitably chosen curve, one shows that the curve cannot contract without its length becoming too large. This yields a nontrivial closed geodesic that is distinct from the shortest one.\\
\hspace*{1cm}In the second case, called \textit{three-legged starfish}, we first assume that all three legs are longer than $\frac{\operatorname{16\sqrt{2}Area}(S^2)}{L_1}$. Using an argument analogous to that in the long sphere case, we obtain a figure‐eight geodesic. This figure‐eight geodesic is supported on two legs that differ from those supporting the shortest closed geodesic. In the alternative case, when one of the legs is short, we construct a sweepout of the sphere with uniformly bounded lengths. Applying the Birkhoff min-max method leads to the existence of a sufficiently short simple closed geodesic.\\
\hspace*{1cm}\textbf{Acknowledgments.} The author would like to express sincere gratitude to the research supervisor, Prof. Yevgeny Liokumovich, for their invaluable guidance and support throughout this research. He is grateful to Prof. Alexander Nabutovsky, Prof. Regina Rotman, Bruno Staffa, Herng Yi Cheng, and Mohammad Alattar for many very valuable discussions. Part of this work was done during the author’s stay in Bonn. He is thankful to the Rheinische Friedrich-Wilhelms-Universität Bonn and the Hausdorff Research Institute for Mathematics for the hospitality. This work was partly supported by the German Research Foundation (DFG).         

\section{General proof} \label{sec-2}
\subsection{The Almgren--Pitts min-max theory} \label{subsec-2-1}
\hspace*{1cm}We begin by introducing the notations from the Almgren--Pitts min-max theory. Let $(M,g)$ be an orientable, compact Riemannian 2-dimensional manifold, possibly with boundary $\partial M$, and assume that $M$ is isometrically embedded in some Euclidean space $\mathbb{R}^L$. Let $X$ be a cubical subcomplex of the $m$-dimensional cube $I^m = [0,1]^m$.\\
\hspace*{1cm} In this section we consider
\begin{itemize}
    \item the space $\mathbf{I}_{1}(M; \mathbb{Z}_2)$ the space of $1$-dimensional mod $2$ flat chains in $\mathbb{R}^L$ whose supports are contained in $M$;
    \item the space $\mathcal{Z}_{1}(M; \mathbb{Z}_2) \bigl(\mathcal{Z}_1(M, \partial M; \mathbb{Z}_2)\bigr)$ of mod $2$ flat chains $T \in \mathbf{I}_{1}(M; \mathbb{Z}_2)$ with $\partial T = 0$ $\bigl(\spt(\partial T) \subset \partial M\bigr)$;
    \item the closure $\mathcal{V}_1(M)$, in the weak topology, of the space of 1-dimensional rectifiable varifolds in $\mathbb{R}^L$ with support contained in $M$.
\end{itemize}
We assume that $\mathbf{I}_1(M; \mathbb{Z}_2), \mathcal{Z}_1(M; \mathbb{Z}_2)$, and $\mathcal{Z}_1(M, \partial M; \mathbb{Z}_2)$ have the
topology induced by the flat metric. When endowed with the topology of the mass norm, these spaces will be denoted by $\mathbf{I}_1(M; \mathbf{M}; \mathbb{Z}_2), \mathcal{Z}_1(M; \mathbf{M}; \mathbb{Z}_2)$, and $\mathcal{Z}_1(M, \partial M; \mathbf{M}; \mathbb{Z}_2)$, respectively. For each $j \in \mathbb{N}$, let $I(1,j)$ denote the cube complex on $I^1$ whose 1-cells and 0-cells (or vertices) are, respectively,
\[
[0, 3 - j],\ [3 - j, 2\cdot 3 - j],\ \ldots,\ [1 - 3 - j, 1] \quad \text{and} \quad [0],\ [3 - j],\ \ldots,\ [1 - 3 - j],\ [1].
\]
Define $I(m,j)$ as the cell complex on $I^m$ by setting
\[
I(m,j) = I(1,j) \otimes \cdots \otimes I(1,j) \quad \text{($m$ times)}.
\]
Then $\alpha = \alpha_1 \otimes \cdots \otimes \alpha_m$ in $I(m,j)$ is a $q$-cell if each $\alpha_i$ is a cell of $I(1,j)$ and
\[
\sum_{i=1}^{m} \dim(\alpha_i) = q.
\]
We will often abuse notation by identifying a $q$-cell $\alpha$ with its support, namely, $\alpha_1 \times \cdots \times \alpha_m \subset I^m$.\\
\hspace*{1cm}The cube complex $X(j)$ is defined as the union of all cells in $I(m,j)$ whose supports lie in some cell of $X$. We denote by $X(j)_q$ the set of all $q$-cells in $X(j)$. Two vertices $x,y \in X(j)_0$ are said to be \textit{adjacent} if they belong to a common 1-cell in $X(j)_1$.\\
\hspace*{1cm}For any $i,j \in \mathbb{N}$, we define the map $\mathbf{n}(i,j) : X(i)_0 \to X(j)_0$ such that $\mathbf{n}(i,j)(x)$ is the vertex in $X(j)_0$ closest to $x$ (see Section 7.1 in~\cite{MarquesNevesWillmore} for a precise definition).\\
\hspace*{1cm}Given a map $\phi : X(j)_0 \rightarrow \mathcal{Z}_1(M;\mathbb{Z}_2)$, its \textit{fineness} is defined by
\[
\mathbf{f}(\phi) = \sup\Bigl\{\mathbf{M}\bigl(\phi(x)-\phi(y)\bigr) : x,y \text{ are adjacent vertices in } X(j)_0\Bigr\}.
\]
This notion provides a discrete measure of continuity with respect to the mass norm.
\begin{definition}
    Let $\phi_i : X(k_i)_0 \rightarrow \mathcal{Z}_1(M; \mathbb{Z}_2), i = 1,2$. We say that $\phi_1$ is \textit{$X$-homotopic to $\phi_2$ in $\mathcal{Z}_1(M; \mathbb{Z}_2)$ with fineness $\delta$} if there exist $k \in \mathbb{N}$ and a map
    \begin{equation*}
        \psi : I(1,k)_0 \times X(k)_0 \rightarrow \mathcal{Z}_1(M; \mathbb{Z}_2)
    \end{equation*}
    satisfying the following conditions:
    \begin{enumerate}
        \item $\mathbf{f}(\psi) < \delta$;
        \item for $i = 1,2$ and any $x \in X(k)_0$,
    \end{enumerate}
    \begin{equation*}
        \psi\bigl([i-1], x\bigr) = \phi_i\bigl(\mathbf{n}(k,k_i)(x)\bigr).
    \end{equation*}
\end{definition}
Instead of considering continuous maps from $X$ into $\mathcal{Z}_1(M; \mathbf{M}; \mathbb{Z}_2)$, Almgren--Pitts theory works with sequences of discrete maps into $\mathcal{Z}_1(M; \mathbb{Z}_2)$ whose fineness tends to zero.
\begin{definition}
    An \textit{$(X, \mathbf{M})$-homotopy sequence of mappings into $\mathcal{Z}_1(M;\mathbf{M};\mathbb{Z}_2)$} is a sequence of mappings $S = \{\phi_i\}_{i \in \mathbb{N}}$,
    \begin{equation*}
        \phi_i : X(k_i)_0 \rightarrow \mathcal{Z}_1(M; \mathbb{Z}_2),
    \end{equation*}
    such that $\phi_i$ is $X$-homotopic to $\phi_{i+1}$ in $\mathcal{Z}_1(M; \mathbf{M}; \mathbb{Z}_2)$ with fineness $\delta_i$, and
    \begin{enumerate}
        \item $\lim_{i \rightarrow \infty} \delta_i = 0$;
        \item $\sup\bigl\{\mathbf{M}\bigl(\phi_i(x)\bigr) : x \in X(k_i)_0, i \in \mathbb{N}\bigr\} < +\infty$.
    \end{enumerate}
\end{definition}
The following definition clarifies what it means for two distinct homotopy sequences of mappings into $\mathcal{Z}_1(M; \mathbf{M}; \mathbb{Z}_2)$ to be homotopic.
\begin{definition}
    Let $S^1 = \{\phi_i^1\}_{i \in \mathbb{N}}$ and $S^2 = \{\phi_i^2\}_{i \in \mathbb{N}}$ be $(X, \mathbf{M})$-homotopy sequences of mappings into $\mathcal{Z}_1(M; \mathbf{M}; \mathbb{Z}_2)$. We say that $S^1$ is \textit{homotopic with} $S^2$ if there exists a sequence $\{\delta_i\}_{i \in \mathbb{N}}$ such that
    \begin{enumerate}
        \item $\phi_i^1$ is $X$-homotopic to $\phi_i^2$ in $\mathcal{Z}_1(M; \mathbf{M}; \mathbb{Z}_2)$ with fineness $\delta_i$;
        \item $\lim_{i \rightarrow \infty} \delta_i = 0$.
    \end{enumerate}
\end{definition}
For a given map $\phi : X_0 \rightarrow \mathcal{Z}_1(M; \mathbb{Z}_2)$, there exists a continuous extension $\Phi: X \rightarrow \mathcal{Z}_1(M; \mathbf{M}; \mathbb{Z}_2)$ that is unique up to homotopy and preserves homotopy classes. The map $\Phi$ is called \textit{the Almgren extension} of $\phi$ (see Theorem 3.10 in~\cite{MarquesNeves} for a precise definition). \\
\hspace*{1cm}The relation "is homotopic with" defines an equivalence relation on the set of all $(X, \mathbf{M})$-homotopy sequences of mappings into $\mathcal{Z}_1(M; \mathbf{M}; \mathbb{Z}_2)$. The equivalence class of any such sequence is called an \textit{$(X, \mathbf{M})$-homotopy class of mappings into $\mathcal{Z}_1(M; \mathbf{M}; \mathbb{Z}_2)$}. The set of all equivalence classes is denoted by $[X,\mathcal{Z}_1(M; \mathbf{M}; \mathbb{Z}_2)]^{\#}$. \\
\hspace*{1cm}Given $\Pi \in [X,\mathcal{Z}_1(M; \mathbf{M}; \mathbb{Z}_2)]^{\#}$, define the function
\begin{equation*}
    \mathbf{L}: \Pi \rightarrow [0, +\infty]
\end{equation*}
by
\begin{equation*}
    \mathbf{L}(S) = \limsup_{i \rightarrow \infty} \max\bigl\{\mathbf{M}\bigl(\phi_i(x)\bigr): x \in \operatorname{dmn}(\phi_i)\bigr\}, \quad \text{where } S = \{\phi_i\}_{i \in \mathbb{N}}.
\end{equation*}
Note that $\mathbf{L}(S)$ serves as the discrete analogue of the maximum length obtained from a continuous map into $\mathcal{Z}_1(M; \mathbf{M}; \mathbb{Z}_2)$.\\
\hspace*{1cm}For a given sequence $S = \{\phi_i\}_{i \in \mathbb{N}} \in \Pi$, consider the compact subset $\mathbf{K}(S)$ of $\mathcal{V}_1(M)$ defined by
\begin{align*}
    \mathbf{K}(S) = \bigl\{V : V = \lim_{j \rightarrow \infty} |\phi_{i_{j}}(x_j)| \text{ as varifolds, for some increasing sequence} \{i_j\}_{j \in \mathbb{N}} \\ \text{ and $x_j \in \operatorname{dmn}(\phi_{i_{j}})$}\bigr\}  
\end{align*}
This set is the discrete replacement to the image of a continuous map into $\mathcal{Z}_1(M; \mathbf{M}; \mathbb{Z}_2)$.
\begin{definition}
    The \textit{width} of $\Pi$ is defined as
    \begin{equation*}
        \mathbf{L}(\Pi) = \inf\{\mathbf{L}(S) : S \in \Pi\}.
    \end{equation*}
    A sequence $S \in \Pi$ is called a \textit{critical sequence} for $\Pi$ if
    \begin{equation*}
        \mathbf{L}(S) = \mathbf{L}(\Pi).
    \end{equation*}
    The \textit{critical set} $\mathbf{C}(S)$ of a critical sequence $S \in \Pi$ is given by
    \begin{equation*}
        \mathbf{C}(S) = \mathbf{K}(S) \cap \bigl\{V : ||V||(M) = \mathbf{L}(S)\bigr\}.
    \end{equation*}
\end{definition}
Now, we introduce the min-max problem in this framework. The Almgren isomorphism theorem provides an isomorphism between 
\begin{equation*}
    \pi_1\bigl(\mathcal{Z}_1(M; \mathbb{Z}_2)\bigr) \quad \text{and} \quad H_{2}(M;\mathbb{Z}_2) = \mathbb{Z}_2.
\end{equation*}
Consequently,
\begin{equation*}
    H^1\bigl(\mathcal{Z}_1(M; \mathbb{Z}_2); \mathbb{Z}_2\bigr) = \mathbb{Z}_2,
\end{equation*}
with generator $\lambda$. We denote by $\lambda^p$ the $p$-fold cup product of $\lambda$ with itself.
\begin{definition}
    Let $\Pi \in [X,\mathcal{Z}_1(M; \mathbf{M}; \mathbb{Z}_2)]^{\#}$. We say that $\Pi$ is a \textit{class of (discrete) $p$-sweepouts} if, for any sequence $S = \{\phi_i\} \in \Pi$, the Almgren extension $\Phi_i: X \rightarrow \mathcal{Z}_1(M; \mathbf{M}; \mathbb{Z}_2)$ of $\phi_i$ is a $p$-sweepout for all sufficiently large $i$.
\end{definition}
We say $X$ is $p$-admissible if there exists a $p$-sweepout $\Phi: X \rightarrow \mathcal{Z}_1(M; \mathbb{Z}_2)$ that has no concentration of mass (Definition 3.7 in~\cite{MarquesNeves}).
\begin{definition}
    Let $\Pi \in [X,\mathcal{Z}_1(M; \mathbf{M}; \mathbb{Z}_2)]^{\#}$. We say that $\Pi$ is a \textit{class of (discrete) $p$-sweepouts} if for any $S = \{\phi_i\} \in \Pi$, the Almgren extension $\Phi_i: X \rightarrow \mathcal{Z}_1(M; \mathbf{M}; \mathbb{Z}_2)$ of $\phi_i$ is a $p$-sweepout for every sufficiently large $i$.
\end{definition}
These definitions are compatible in the sense that a continuous map $\Phi$ is a $p$-sweepout if and only if its discretization $\Pi$ is a class of $p$-sweepouts (see Lemma 4.6 in~\cite{MarquesNeves}).
\begin{definition}
    Let $\mathcal{D}_p$ be the collection of all classes of $p$-sweepouts,
    \begin{equation*}
        \Pi \in [X,\mathcal{Z}_1(M; \mathbf{M}; \mathbb{Z}_2)]^{\#},
    \end{equation*}
    where $X$ is any $p$-admissible cubical subcomplex. The \textit{$p$-width} of $M$ is defined by
    \begin{equation*}
        \omega_p(M) = \inf_{\Pi \in \mathcal{D}_p} \mathbf{L}(\Pi).
    \end{equation*}
\end{definition}
It is not evident a priori whether $\omega_p(M)$ coincides with the width $\mathbf{L}(\Pi)$ for some class of $p$-sweepouts $\Pi$. The following proposition addresses the situation in which this equality fails.
\begin{proposition} \label{widthless}
    Suppose there exists $p \in \mathbb{N}$ such that, for every $p$-admissible $X$, we have
    \begin{equation*}
        \omega_p(M) < \mathbf{L}(\Pi) \quad \text{for every class of $p$-sweepouts } \Pi \in [X,\mathcal{Z}_1(M; \mathbf{M}; \mathbb{Z}_2)]^{\#}.
    \end{equation*}
    Then, for any $\varepsilon > 0$, there exist infinitely many immersed closed geodesics in $M$ with lengths bounded above by $\omega_p(M) + \varepsilon$.
\end{proposition}
\textbf{Proof.} By Lemma 4.7 in~\cite{MarquesNeves}, one may find sequences of $p$-admissible cubical subcomplexes $X_k$ and classes of $p$-sweepouts $\Pi_k \in [X_k,\mathcal{Z}_1(M; \mathbf{M}; \mathbb{Z}_2)]^{\#}$ such that
\begin{equation*}
    \mathbf{L}(\Pi_1) > \mathbf{L}(\Pi_2) > \cdots > \mathbf{L}(\Pi_k) > \mathbf{L}(\Pi_{k+1}) > \cdots,
\end{equation*}
and
\begin{equation*}
    \lim_{k \rightarrow \infty} \mathbf{L}(\Pi_k) = \omega_p(M).
\end{equation*}
Then, by Propositions 2.12, 2.13, and 4.1 in~\cite{ChodoshMantoulidis}, there exist closed geodesics $\sigma_{k,1}, \ldots, \sigma_{k,N_k}$ (allowing repetitions) such that
\begin{equation*}
    \mathbf{L}(\Pi_k) = \sum_{j = 1}^{N_k} L(\sigma_{k,j}).
\end{equation*}
Since $N_k$ is uniformly bounded from above (because $L(\sigma_{k,j}) \geq 2 \operatorname{inj}(M)$ for each $j$), one may extract a subsequence for which the $\sigma_{k,j}$ converge, yielding the desired geodesics. \hfill\qedsymbol{} \\
\hspace*{1cm}The following theorem is an adaptation of Theorem 6.1 in~\cite{MarquesNeves} for the case of 2-dimensional manifolds based on the Lusternik--Schnirelmann theory.
\begin{theorem} \label{widthsequal}
    If $\omega_p(M) = \omega_{p+1}(M)$ for some $p \in \mathbb{N}$, then for any $\varepsilon > 0$, there exist infinitely many immersed closed geodesics in $M$ with lengths bounded above by $\omega_p(M) + \varepsilon$.
\end{theorem}
\textbf{Proof.} By Proposition~\ref{widthless}, we may assume that there exists a $(p+1)$-admissible cubical subcomplex $X$ and a class of $(p+1)$-sweepouts
\begin{equation*}
    \Pi \in [X,\mathcal{Z}_1(M; \mathbf{M}; \mathbb{Z}_2)]^{\#}
\end{equation*}
such that $\omega_{p+1}(M) = \mathbf{L}(\Pi)$. Then, by Proposition 2.7 in~\cite{MarquesNeves}, one can find a critical sequence $S = \{\phi_i\}_{i \in \mathbb{N}} \in \Pi$ for which every $\Sigma \in \mathbf{C}(S)$ is a stationary varifold with mass
\begin{equation*}
    \mathbf{L}(S) = \mathbf{L}(\Pi) = \omega_{p+1}(M).
\end{equation*}
If $\Phi_i: X \rightarrow \mathcal{Z}_1(M; \mathbf{M}; \mathbb{Z}_2)$ denotes the Almgren extension of $\phi_i$, then the fact that $\Pi$ is a class of $(p+1)$-sweepouts implies that $\Phi_i$ is a $(p+1)$-sweepout for all sufficiently large $i$.\\
\hspace*{1cm}Let $\mathcal{S}$ be the set of all stationary integral varifolds with length not exceeding $\omega_{p+1}(M)$ whose support is a union of immersed closed geodesics. Similarly, let $\mathcal{T}$ denote the set of all mod $2$ flat chains $T \in \mathcal{Z}_1(M; \mathbb{Z}_2)$ with $\mathbf{M}(T) \leq \omega_{p+1}(M)$ such that either $T = 0$ or the support of $T$ is a union of immersed closed geodesics. Suppose, for contradiction, that both $\mathcal{S}$ and $\mathcal{T}$ are finite. The remainder of the proof follows the same structure as that of Theorem 6.1 in~\cite{MarquesNeves}; we now outline the main ideas. \\
\hspace*{1cm}Finiteness of $\mathcal{S}$ and $\mathcal{T}$ allows us to construct cubical subcomplexes $Y_i$ such that, for some $\eta > 0$,
\begin{equation*}
    \mathbf{F}(|\phi_i(x)|, \mathcal{S}) \geq \eta,
\end{equation*}
and
\begin{equation*}
    \mathbf{F}(|\Phi_i(x)|, \mathcal{S}) < 2\eta \quad \text{for all } x \in X \setminus Y_i,
\end{equation*}
for sufficiently large $i$, where $\mathbf{F}$ denotes the $\mathbf{F}$-metric on $\mathbf{I}_1(M;\mathbb{Z}_2)$, and the restriction $(\Phi_i)_{|Y_i}$ is a $p$-sweepout. Consider the sequence $\tilde{S} = \{\psi_i\}$ defined by
\begin{equation*}
    \psi_i = (\phi_i)_{|Y_i}: (Y_i)_0 \rightarrow \mathcal{Z}_1(M; \mathbb{Z}_2),
\end{equation*}
and set
\begin{equation*}
    L = \mathbf{L}(\tilde{S}) = \limsup_{i \rightarrow \infty} \max\bigl\{\mathbf{M}\bigl(\psi_i(y)\bigr): y \in (Y_i)_0\bigr\}.
\end{equation*}
Clearly, $L \leq \omega_{p+1}(M) = \omega_p(M)$. There are two possibilities: either $L < \omega_p(M)$ or $L = \omega_p(M)$. In the former case, since $\Phi_i$ is a $p$-sweepout for sufficiently large $i$, this contradicts the definition of $\omega_p(M)$. In the latter case, the inclusion
\begin{equation*}
    \mathbf{C}(\tilde{S}) \subset \{V: \mathbf{F}(V, \mathcal{S}) \geq \eta\}
\end{equation*}
implies that, although every element of $\mathbf{C}(\tilde{S})$ is stationary, none is almost minimizing in annuli. Consequently, there exists a sequence $\tilde{S}^* = \{\psi_i^*\}$ of maps
\begin{equation*}
    \psi_i^*: Y_i(l_i)_0 \rightarrow \mathcal{Z}_1(M; \mathbb{Z}_2)
\end{equation*}
such that:
\begin{itemize}
    \item $\psi_i$ and $\psi^*_i$ are $X$-homotopic with fineness tending to zero as $i \rightarrow \infty$,
    \item $\mathbf{L}(\tilde{S}^*) < \mathbf{L}(\tilde{S}) = L = \omega_p(M)$.
\end{itemize}
Since $\{\Psi_i^*\}$ also forms a $p$-sweepout, the inequality $\mathbf{L}(\tilde{S}^*) < \omega_p(M)$ contradicts the definition of $\omega_p(M)$. Thus, both cases lead to a contradiction, and it follows that there must be infinitely many distinct immersed closed geodesics in $M$.\hfill\qedsymbol{}

\subsection{Proof of Main Theorem} \label{subsec-2-2}
\hspace*{1cm}Let $(S^2, g)$ be a Riemannian 2-dimensional sphere. Denote $A = \Area(S^2, g)$ and $D = \diam(S^2, g)$.
\begin{theorem}
    For any Riemannian metric $g$ on a 2-dimensional sphere $S^2$, there exist two distinct closed geodesics with lengths $L_{1}$ and $L_{2}$ satisfying
    \begin{equation*}
        L_{1} L_{2} \leq 2^9 \cdot 10^4\, A.
    \end{equation*}
\end{theorem}
\textbf{Proof.} Consider the $p$-widths $\omega_{p}$ of $(S^2,g)$ for $p \in \mathbb{Z}_{+}$. By Theorem~\ref{widthsequal}, we can assume that the sequence is strictly increasing. According to Theorem 1.2 in~\cite{ChodoshMantoulidis}, the $\omega_{p}$ correspond to the lengths of unions of immersed closed geodesics on $(S^2, g)$ (possibly with multiplicities). Moreover, by Theorem 1.1 in~\cite{Liokumovich},
\begin{equation} \label{ineq-1}
    \omega_{p} \leq 1600 \sqrt{p A}.
\end{equation}
Denote by $\gamma_{1}$ the shortest closed geodesic and by $L_{1}$ its length. By Theorem 1 in~\cite{Rotman},
\begin{equation*}
    L_{1} \leq 4 \sqrt{2A}.
\end{equation*}
Then, there exists $n \in \mathbb{Z}_+$ such that
\begin{equation} \label{ineq-2}
    \frac{4\sqrt{2A}}{n+1} < L_{1} \leq \frac{4\sqrt{2A}}{n}.
\end{equation}
Suppose that $\omega_{p}$'s correspond to coverings of $\gamma_{1}$. This means that
\begin{equation*}
    \omega_{p} = \varphi(p) L_{1} \quad \text{with} \quad \varphi(p) \geq p.
\end{equation*}
Inequalities~\eqref{ineq-1} and~\eqref{ineq-2} imply that
\begin{equation*}
    \frac{4p\sqrt{2A}}{n+1} < pL_{1} \leq \varphi(p) L_{1} = \omega_{p} \leq 1600 \sqrt{p A},
\end{equation*}
which is equivalent to
\begin{equation*}
    \sqrt{p} < 200\sqrt{2}\,(n+1).
\end{equation*}
Thus, for 
\[
p = \bigl(200\sqrt{2}\,(n+1)\bigr)^2,
\]
the $p$-width does not correspond to a covering of $\gamma_{1}$. Thus, we obtain a closed geodesic $\gamma_{2}$ with length $L_{2} = \omega_{p}$ that is distinct from $\gamma_{1}$. Consequently,
\begin{equation*}
    L_{1} L_{2} \leq \frac{4\sqrt{2A}}{n} \cdot 1600\sqrt{p} \sqrt{A} = \frac{6400\sqrt{2}}{n} \cdot 200\sqrt{2}\,(n+1) A \leq 2^9 \cdot 10^4\, A.
\end{equation*}
The theorem is proved.\hfill\qedsymbol{}
\section{Better bounds} \label{sec-3}
\subsection{Long sphere} \label{subsec-3-1}
Let the shortest closed geodesic $\gamma_1$ on $S^2$ be simple and have Morse index one. By the Jordan separation theorem, $\gamma_1$ divides $S^2$ into two regions. Choose points $x$ and $y$ in each of them on the maximal distance from $\gamma_1$. Denote the corresponding regions by $M_x$ and $M_y$, and let $d(x, \gamma_1) = d_x$ and $d(y, \gamma_1) = d_y$.
\begin{theorem}
    Suppose $d_x, d_y > \frac{170A}{L_{1}}$. Then there exists a closed geodesic $\gamma_{2}$ distinct from $\gamma_{1}$ with length $L_{2}$ such that 
    \begin{equation*}
        L_{1} L_{2} \leq 320A.
    \end{equation*}
\end{theorem}
\textbf{Proof.}
    Let $\tau \: [0, D] \rightarrow S^2$ be a minimizing geodesic connecting $x$ and $y$ parameterized with arclength. Now $B(x, d_{x}) \cap B(y, d_{y})$ has measure zero, so by the Coarea inequality
    \begin{equation*}
        A \geq \Area\bigl(B(x, d_{x})\bigr) + \Area\bigl(B(y, d_{y})\bigr) \geq \int_{0}^{\frac{2A}{L_{1}}} L\bigl(S(x, d_{x} - u)\bigr) du + \int_{0}^{\frac{2A}{L_{1}}} L\bigl(S(y, d_{y} - u)\bigr) du.
    \end{equation*}
    Note that
    \begin{equation*}
        \int_{0}^{\frac{2A}{L_{1}}} \frac{L_{1}}{2} du = A.
    \end{equation*}
    Hence, there is $u \in [0, \frac{2A}{L_{1}}]$ such that
    \begin{equation*}
        L\bigl(S(x, d_{x} - u)\bigr) + L\bigl(S(y, d_{y} - u)\bigr) \leq \frac{L_{1}}{2}
    \end{equation*}
    and $S(x, d_{x} - u)$ and $S(y, d_{y} - u)$ are disjoint unions of simple closed curves. Denote $\alpha_{0}$ the component of $S(x, d_{x} - u)$ through $\tau(d_{x} - u)$ and $\alpha_{1}$ the component of $S(y,d_{y} - u)$ through $\tau(D - d_{y} + u)$. \\
    \hspace*{1cm}Since $\gamma_{1}$ is simple, the Jordan separation theorem implies that it divides $S^2$ into two regions, $M_{x}$ and $M_{y}$, with $x \in M_{x}$ and $y \in M_{y}$. Moreover, as $\gamma_{1}$ has Morse index 1, for any $\varepsilon > 0$ there exist closed curves $\gamma_{x}^{\varepsilon}$ and $\gamma_{y}^{\varepsilon}$ satisfying:
    \begin{itemize}
        \item $d(\gamma_{1}, \gamma_{x}^{\varepsilon}) < \varepsilon$ and $d(\gamma_{1}, \gamma_{y}^{\varepsilon}) < \varepsilon$,
        \item $L\bigl(\gamma_{x}^{\varepsilon}\bigr) < L_{1}$ and $L\bigl(\gamma_{y}^{\varepsilon}\bigr) < L_{1}$,
        \item $\gamma_{x}^{\varepsilon} \subset M_{x}$ and $\gamma_{y}^{\varepsilon} \subset M_{y}$.
    \end{itemize}
    Then, if we apply the Birkhoff curve shortening process to $\gamma_{x}^{\varepsilon}$, we will obtain a homotopy between $\gamma_{x}^{\varepsilon}$ and a point curve or a closed geodesic different from $\gamma_{1}$. In the latter case, we are done. Let's assume the former case. Denote by $\sigma_{0}$ the curve in the homotopy which intersects $\alpha_{0}$ the last. It has length less than $L(\alpha_{0}) \leq \frac{L_{1}}{2}$. Also, by the triangle inequality, we have $d(x, \sigma_{0}) > d(x, \gamma_{1}) - \frac{2A}{L_{1}} - \frac{L(\alpha_{0})}{2} > \frac{170A}{L_{1}} - \frac{2A}{L_{1}} - \frac{8A}{L_{1}} = \frac{160A}{L_{1}}$. We can repeat the same argument for $\gamma_{y}^{\varepsilon}$ and obtain a homotopy between $\gamma_{y}^{\varepsilon}$ and curve $\sigma_{1}$ with the same properties. Taking $\varepsilon \rightarrow 0$ and concatenating these two homotopies, we obtain a homotopy $h_{s}$ such that $h_{0} = \sigma_{0}$ and $h_{1} = \sigma_{1}$. \\
    \begin{figure}
        \centering
        \includegraphics[scale=0.5]{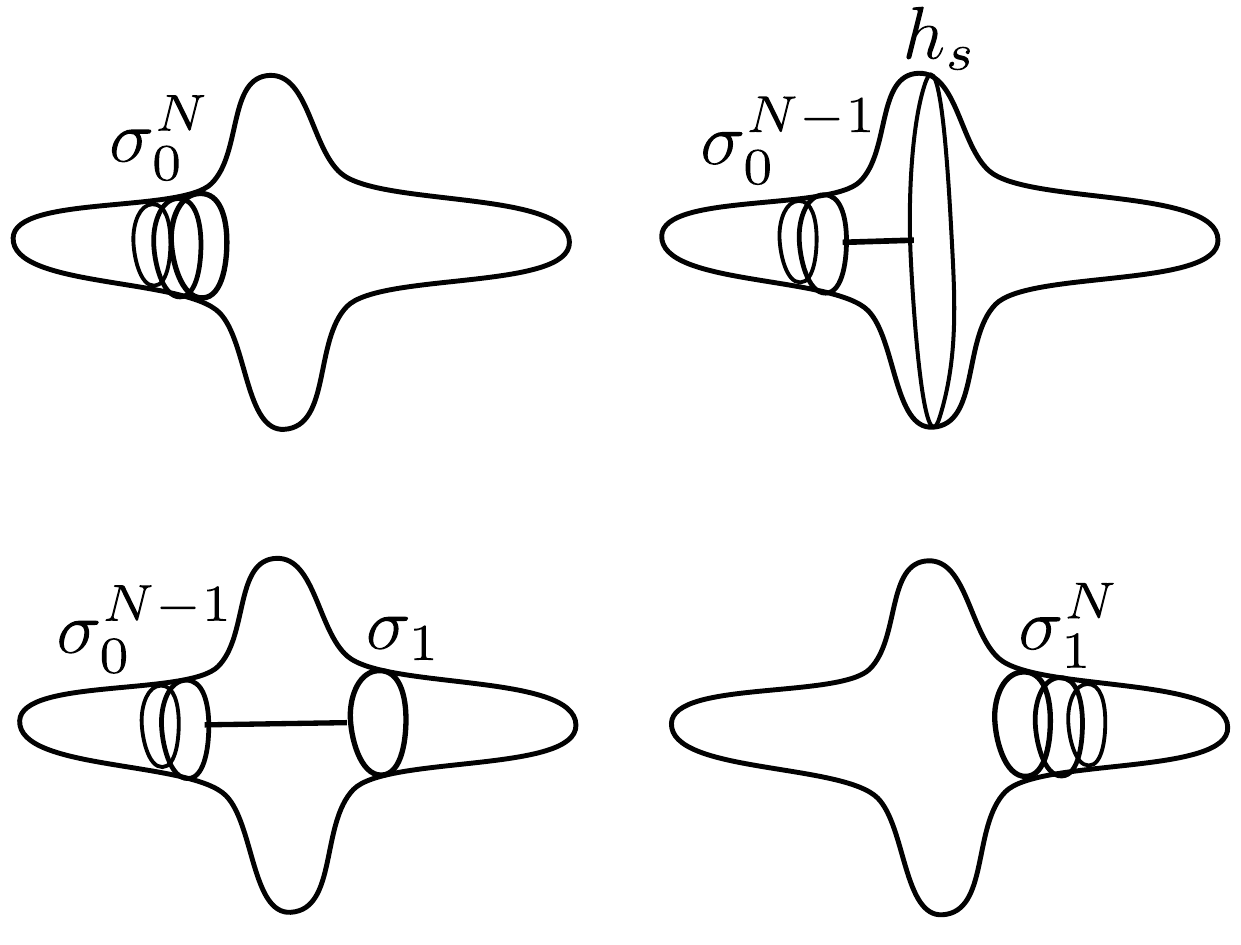}
        \caption{Homotopy $H_{t}^{N}$}
        \label{H_t^N}
    \end{figure}
    \hspace*{1cm}Now we want to describe homotopy $H_{t}^N$ between $\sigma_{0}^N$ and $\sigma_{1}^N$ where $\sigma_{i}^N$ denotes the $N$-fold covering of $\sigma_{i}$ for some integer $N$. We pull each loop $\sigma_{0}$ one by one to $\sigma_{1}$ using the homotopy $h_{s}$ (see Fig.~\ref{H_t^N}). During the deformation, the loops are connected through $\tau$. We want to choose $N$ such that
    \begin{equation} \label{ineq-H_t^n}
        L(H_{t}^N) < 2N L_{1}
    \end{equation}
    for each $t$. There exists $n \in \mathbb{N}^*$ such that
    \begin{equation*}
        \frac{4\sqrt{2A}}{n+1} < L_{1} \leq \frac{4\sqrt{2A}}{n}.
    \end{equation*} Notice that the segment of $\tau$ has the length less than or equal to
    \begin{equation*}
        \frac{4A}{L_{1}} + \frac{L_{1}}{2} \leq \bigl(\frac{(n + 1)^2}{8} + \frac{1}{2}\bigr)L_{1}.
    \end{equation*}
    Then
    \begin{align*}
        L(H_{t}^N) \leq (2N - 1) \frac{L_{1}}{2} + L_{1} + \bigl(\frac{(n + 1)^2}{8} + \frac{1}{2}\bigr)L_{1} = 2NL_{1} - (N - \frac{1}{2}) L_{1} + \bigl(\frac{(n + 1)^2}{8} + \frac{1}{2}\bigr)L_{1}.
    \end{align*}
    Then $L(H_{t}^N) < 2N L_{1}$ holds if
    \begin{equation*}
        N > \frac{(n+1)^2}{8} + 1.
    \end{equation*}
    So, we pick $N = \floor{{\frac{(n+1)^2}{8} + 2}}$. Denote by $C_0$ the interior of the region bounded by $\sigma_0$ that contains $x$, and similarly define $C_1$ for $\sigma_1$ and $y$.  Since each $\sigma_i$ is obtained via the Birkhoff curve shortening process, $C_i$ is convex for $i\in\{0,1\}$.  Note that $C_0\cap C_1=\emptyset$.  Let $\mathcal{D}$ denote the Birkhoff curve shortening map.  For any integer $k\ge1$, set
    \begin{equation*}
        \mathcal{D}^k\gamma := \mathcal{D}\bigl(\mathcal{D}^{k-1}\gamma\bigr),
    \end{equation*}
    with $\mathcal{D}^1\gamma=\mathcal{D}\gamma$.
    For $i \in \{0,1\}$, define
    \begin{equation*}
        I_{i} = \Bigl\{\,t \in [0,1]\Big | \exists\,k\in\mathbb{N}:\,\mathcal{D}^k\bigl(H_{t}^N\bigr)\subset C_{i}\Bigr\}.
    \end{equation*}
    We claim that each $I_i$ is open.  Indeed, fix $i \in \{0,1\}$ and choose $t\in I_{i}$.  By definition, there exists $k\in\mathbb{N}$ such that 
    \begin{equation*}
        \mathcal{D}^k\bigl(H_{t}^N\bigr) \subset C_{i}.
    \end{equation*}
    Since $C_{i}$ is open, we can find $\varepsilon>0$ with
    \begin{equation*}
        U_{\varepsilon}\bigl(\mathcal{D}^k(H_{t}^N)\bigr) \subset C_{i},
    \end{equation*}
    where $U_{\varepsilon}(\eta)$ denotes the $\varepsilon$-neighborhood of a curve $\eta$.  By the continuity of the Birkhoff curve shortening map (see chapter 5 in~\cite{ColdingMinicozzi}), there exists $\delta>0$ such that whenever two closed curves $\eta_{1},\eta_{2}$ satisfy
    \begin{equation*}
        d\bigl(\eta_{1},\eta_{2}\bigr) < \delta,
    \end{equation*}
    then
    \begin{equation*}
        d\bigl(\mathcal{D}^k(\eta_{1}),\,\mathcal{D}^k(\eta_{2})\bigr) < \varepsilon.
    \end{equation*}
    By continuity of $t \mapsto H_{t}^N$, there is $h > 0$ so that for all $t'\in [\,t-h,\,t+h\,]$,
    \begin{equation*}
        d\bigl(H_{t}^N,\,H_{t'}^N\bigr) < \delta.
    \end{equation*}
    Hence for each $t'\in [\,t-h,\,t+h\,]$,
    \begin{equation*}
        d\bigl(\mathcal{D}^k(H_{t}^N),\,\mathcal{D}^k(H_{t'}^N)\bigr) < \varepsilon,
    \end{equation*}
    which implies 
    \begin{equation*}
        \mathcal{D}^k\bigl(H_{t'}^N\bigr)\subset U_{\varepsilon}\bigl(\mathcal{D}^k(H_{t}^N)\bigr)
        \subset C_{i}.
    \end{equation*}
    Therefore $t'\in I_{i}$ for all $t'$ in a neighborhood of $t$, proving that $I_{i}$ is open. Since $C_{0} \cap C_{1} = \emptyset$ and $\mathcal{D}$ respects convex sets, we have $I_{0} \cap I_{1} = \emptyset$. Since $0 \in I_{0}$ and $1 \in I_{1}$, there exists $t \in [0,1] - (I_{0} \cup I_{1})$. Then $\mathcal{D}^k(H_{t}^N) \cap \bigr(S^2 - (C_{0} \cup C_{1})\bigl) \neq \emptyset$ for each $k$ (cf.~\cite{Bangert}). Denote $\sigma = H_{t}^N$. \\
    \hspace*{1cm}Applying the Birkhoff curve shortening process to $\sigma$, we either obtain a closed geodesic or a point curve. Notice that during the process, it cannot pass through either $x$ or $y$ because
    \begin{align*}
        L(\sigma) < 2NL_{1} \leq 2\bigl(\frac{(n+1)^2}{8} + 2\bigr) \frac{4\sqrt{2}\sqrt{A}}{n} \leq \bigl(\frac{(n+1)^2}{8} + 2\bigr) \frac{128A}{n^2L_{1}} = \\ = (16 + \frac{32}{n} + \frac{16}{n^2} + \frac{256}{n^2}) \frac{A}{L_{1}} \leq \frac{320A}{L_{1}}
    \end{align*}
    which contradicts the inequality
    \begin{equation*}
        2 d\bigl(x, S^2 - (C_{0} \cup C_{1})\bigr) = 2 d(x, \sigma_{0}) > \frac{320A}{L_{1}}.
    \end{equation*}
    So, $\sigma$ cannot be contracted to a point. Also, it cannot converge to a covering of $\gamma_{1}$ due to (\ref{ineq-H_t^n}). Denote the obtained closed geodesic by $\gamma_{2}$ distinct from $\gamma_1$ and its length by $L_{2}$. Then
    \begin{equation*}
        L_{1} L_{2} < 320A.
    \end{equation*}
    The theorem is proved. \hfill\qedsymbol{}
\subsection{Three-legged starfish} \label{subsec-3-2}
\hspace*{1cm}Let's discuss the case of the shortest closed geodesic $\gamma_{1}$ being stable with one self-intersection. In this case $\gamma_{1}$ divides $S^2$ into 3 regions. We pick points $x, y$ and $z$ in each of them on the maximal distance from $\gamma_{1}$. Denote corresponding regions by $M_{x}, M_{y}$ and $M_{z}$, where $\overline{M}_{x} \cap \overline{M}_{y}$ is a point. Also, denote $\gamma_{1} \cap \overline{M}_{x} = \gamma_{x}$ and $\gamma_{1} \cap \overline{M}_{y} = \gamma_{y}$. Let $d_{x} = d(x,\gamma_{1})$ (similarly for $d_{y}$ and $d_{z}$).
\begin{theorem}
    Suppose $d_{x}, d_{y}, d_{z} \geq \frac{16\sqrt{2} A}{L_{1}}$. Then there exists a closed geodesic $\gamma_{2}$ distinct from $\gamma_{1}$ with length $L_{2}$ such that 
    \begin{equation*}
        L_{1} L_{2} \leq 16\sqrt{2}A.
    \end{equation*}
\end{theorem}
\textbf{Proof.}
    Denote $\tau_{y}^x, \tau_{z}^y$ and $\tau_{x}^z$ minimizing geodesics between corresponding points. Let $\tau \: [0, L] \rightarrow S^{2}$ be $\tau_{x}^{z}$ with the arclength parameter $t$. Notice that $\frac{L_{1}}{8\sqrt{2}} \leq \frac{2\sqrt{2}A}{L_{1}} < d_{x}$ and $d_{z}$. Also, $L > d_{x} + d_{z}$. Since $B(x, d_{x}) \cap B(z, d_{z})$ has measure zero so by the Coarea inequality
    \begin{equation*}
        A \geq \Area\bigl(B(x, d_{x})\bigr) + \Area\bigl(B(z, d_{z})\bigr) \geq \int_{0}^{\frac{L_{1}}{8\sqrt{2}}} L\bigl(S(x, d_{x} - s)\bigr)ds + \int_{0}^{\frac{L_{1}}{8\sqrt{2}}} L\bigl(S(z, d_{z} - s)\bigr)ds.
    \end{equation*}
    Note that
    \begin{equation*}
        \int_{0}^{\frac{L_{1}}{8\sqrt{2}}} \frac{8\sqrt{2}A}{L_{1}} ds = A.
    \end{equation*}
    Hence, there is a generic $s \in [0, L_{1}]$ such that
    \begin{figure}
        \centering
        \includegraphics[scale=0.45]{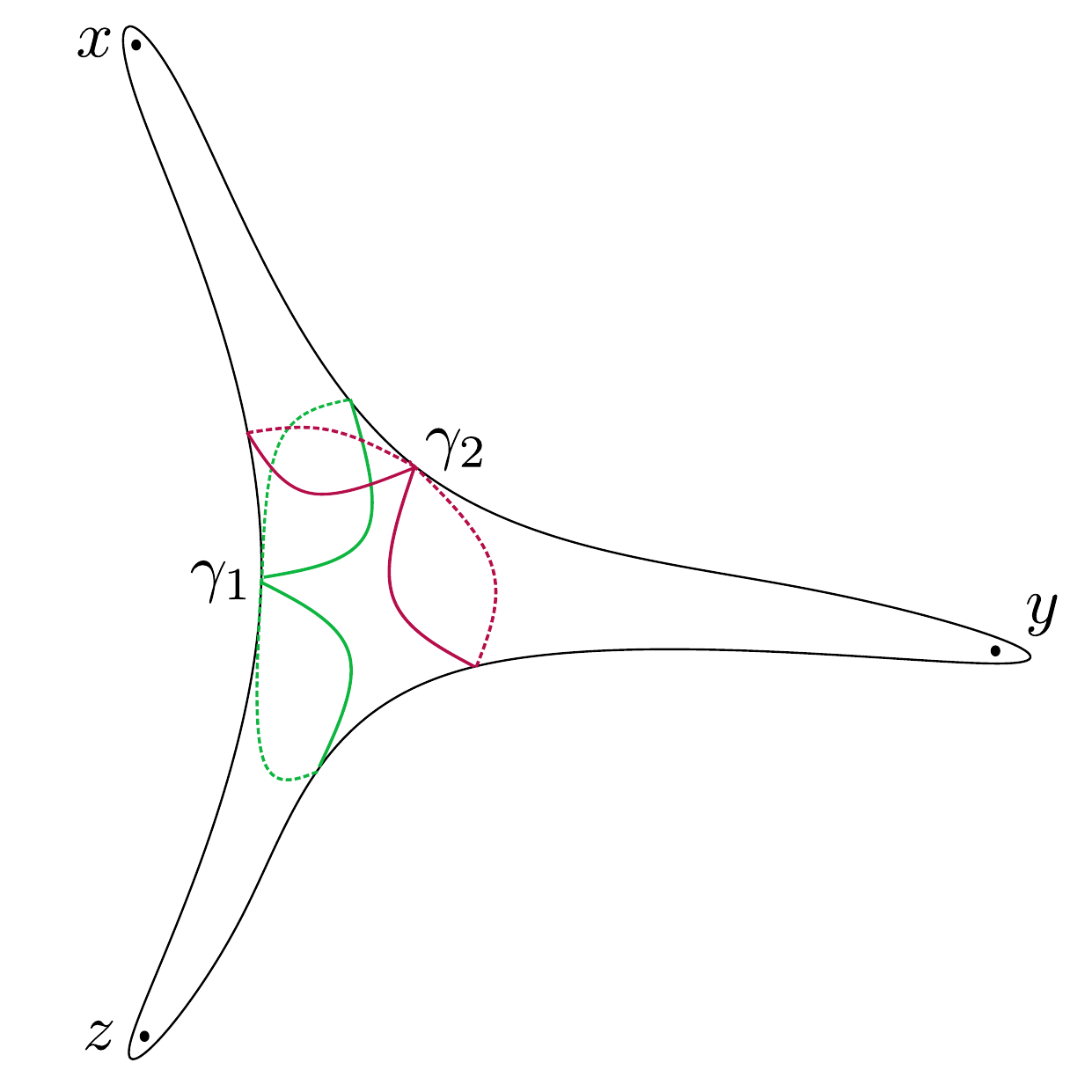}
        \caption{Case of a "long three-legged starfish"}
        \label{long_starfish}
    \end{figure}
    \begin{equation*}
        L\bigl(S(x, d_{x} - s)\bigr) + L\bigl(S(z, d_{z} - s)\bigr) \leq \frac{8\sqrt{2}A}{L_{1}}
    \end{equation*}
    and both $S(x, d_{x} - s)$ and $S(z, d_{z} - s)$ are disjoint unions of simple closed curves. Also, we can assume that $S(x, d_{x} - s)$ intersects $\tau_{y}^x, \tau_{x}^z$ transversely and $S(z, d_{z} - s)$ intersects $\tau_{x}^z, \tau_{z}^y$ transversely. Let $\sigma_{1}$ be the component of $S(x, d_{x} - s)$ through $\tau(d_{x} - s)$ and $\sigma_{2}$ the component of $S(z, d_{z} - s)$ through $\tau(L - d_{z} + s)$. Notice also that $\sigma_{1} \cap \tau_{z}^{y} = \emptyset$ and $\sigma_{2} \cap \tau_{y}^{x} = \emptyset$. By the Jordan Curve Theorem both $\sigma_{1}$ and $\sigma_{2}$ separate $S^2$ into $2$ regions. In the case of $\sigma_{1}$ there must be $x$ on one side and $y, z$ on the other. While in the case of $\sigma_{2}$ there must be $z$ on one side and $x, y$ on the other. Define
    \begin{equation*}
        \sigma = \sigma_{1} \cup \tau |_{[d_{x} - s, L - d_{z} + s]} \cup \sigma_{2} \cup -\tau |_{[d_{x} - s, L - d_{z} + s]}.
    \end{equation*}
    We see that $L(\sigma) \leq \frac{8\sqrt{2}A}{L_{1}} + 4s \leq \frac{8\sqrt{2}A}{L_{1}} + 4\frac{L_{1}}{8\sqrt{2}} \leq \frac{8\sqrt{2}A}{L_{1}} + \frac{8\sqrt{2}A}{L_{1}} = \frac{16\sqrt{2} A}{L_{1}}$. We choose orientation on $\sigma_{1}$ and $\sigma_{2}$ such that the oriented intersection number of $\sigma$ with $\tau$ equals to $2$. So, $\sigma$ is nontrivial. Applying the Birkhoff Curve Shortening Process to $\sigma$ we obtain either a closed geodesic of length $\leq L(\sigma) \leq \frac{16\sqrt{2} A}{L_{1}}$ or a point curve. But it cannot lead to a point curve for if it did then by intersection number argument some curve $\sigma_{s_{0}}$ in the homotopy would have passed through a vertex $x, y$ or $z$ while still intersecting the opposite geodesic $\tau_{y}^x, \tau_{z}^y$ or $\tau_{x}^z$ respectively. But it contradicts with the fact that $L(\sigma_{s_{0}}) \leq \frac{16\sqrt{2} A}{L_{1}}$ and $d(z, \tau_{y}^x), d(x, \tau_{z}^y)$ and $d(y, \tau_{x}^z) > \frac{16\sqrt{2} A}{L_{1}}$. Denote the obtained closed geodesic by $\gamma_{2}$ (see Fig.~\ref{long_starfish}). \\
    \hspace*{1cm}Note that we can choose orientation on $\gamma_{1}$ such that the mod $2$ winding number of $\gamma_{1}$ equals to $1$ around $z$ and to $0$ around $y$. While the mod 2 winding number of $\gamma_{2}$ equals to $1$ around $y$ and to $0$ around $z$. Then $\gamma_{1}$ and $\gamma_{2}$ are distinct and
    \begin{figure}
        \centering
        \includegraphics[scale=0.45]{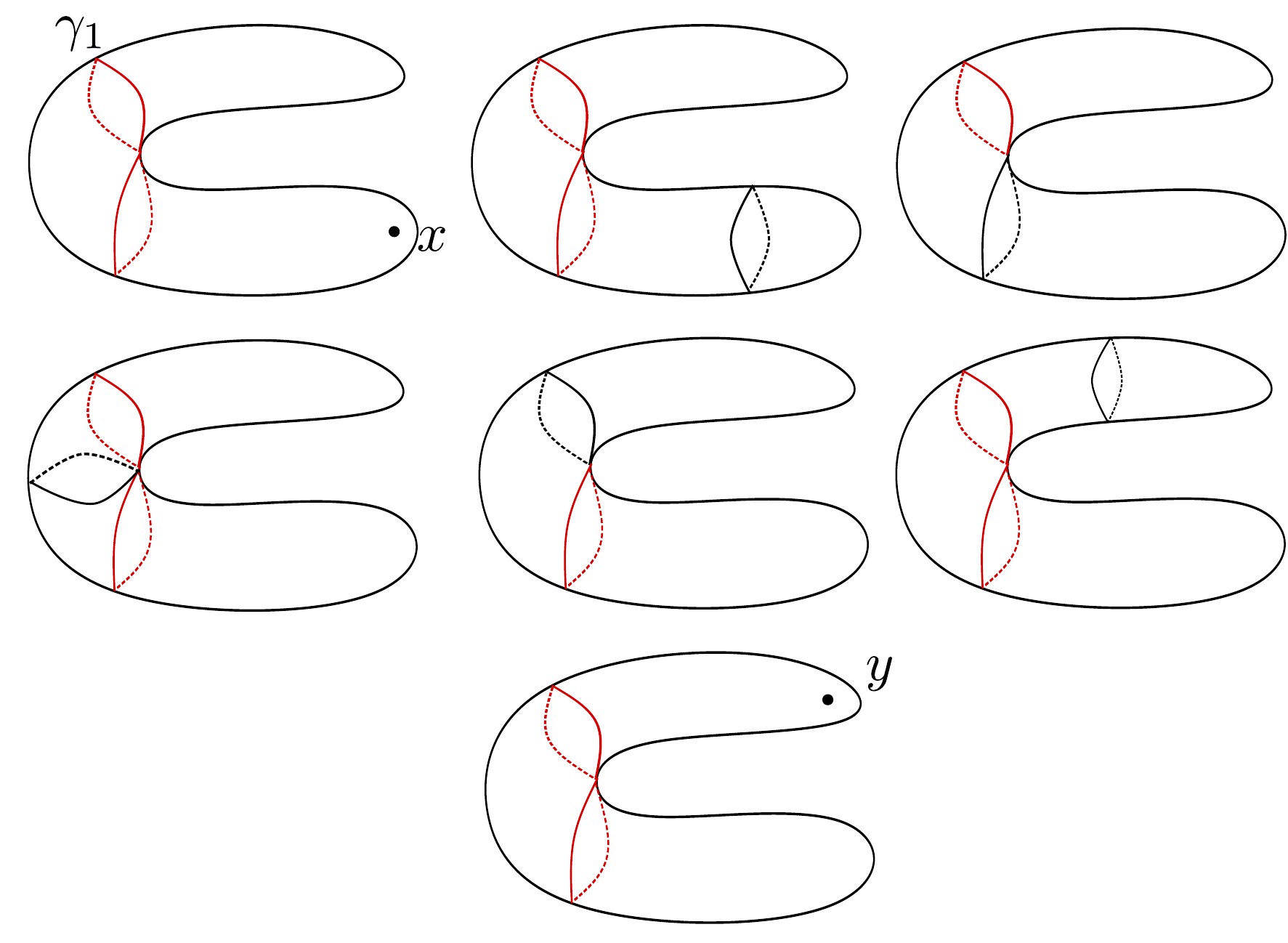}
        \caption{Sweep-out $f$}
        \label{short_starfish1}
    \end{figure}
    \begin{equation*}
        L_{1} L_{2} \leq L_{1} \frac{16\sqrt{2} A}{L_{1}} = 16\sqrt{2} A.
    \end{equation*}
    The theorem is proved. \hfill\qedsymbol{}
\begin{theorem}
    Suppose $d_{z} < \frac{16\sqrt{2} A}{L_{1}}$. Then there exists a simple closed geodesic $\gamma_{2}$ distinct from $\gamma_{1}$ with length $L_{2}$ such that 
    \begin{equation*}
        L_{1} L_{2} \leq (2804\sqrt{2} + 64) A.
    \end{equation*}
\end{theorem}
\textbf{Proof.}
    We will construct a non-contractible map $f \: S^1 \rightarrow \Lambda_{L} S^2$, where $\Lambda_{L} S^2$ is the space of closed continuous curves of length bounded from above by $L$. Using the Birkhoff min-max method (see~\cite{ChambLio}), it will prove existence of a simple closed geodesic with length less or equal than $L$. \\
    \hspace*{1cm}Notice that for any $\varepsilon > 0$, a small perturbation of $\gamma_{1}$ produces closed curves $\gamma_{x}^{\varepsilon}$ and $\gamma_{y}^{\varepsilon}$ satisfying:
    \begin{itemize}
        \item $\gamma_{x}^{\varepsilon} \subset M_{x}$ and $\gamma_{y}^{\varepsilon} \subset M_{y}$,
        \item $L\bigl(\gamma_{x}^{\varepsilon}\bigr)$ and $L\bigl(\gamma_{y}^{\varepsilon}\bigr) \leq L_{1}$,
        \item $d\bigl(\gamma_{x}^{\varepsilon},\gamma_{1}\bigr)$ and $d\bigl(\gamma_{y}^{\varepsilon},\gamma_{1}\bigr) < \varepsilon$.
    \end{itemize}
    Applying the Birkhoff curve shortening process to $\gamma_{x}^{\varepsilon}$ will give us either a closed geodesic different from $\gamma_{1}$ or a homotopy of $\gamma_{x}^{\varepsilon}$ with a point curve inside $M_{x}$. In the former case we are done. In the latter case we can assume the point curve is $x$. Taking $\varepsilon \rightarrow 0$ we obtain a homotopy of $\gamma_{x}$ with the point curve $x$ inside $M_{x}$. Repeating the same argument for $\gamma_{y}^{\varepsilon}$ we obtain a homotopy of $\gamma_{y}$ with the point curve $y$ inside $M_{y}$. Notice that the lengths of the curves in the homotopies don't exceed $L_{1}$.\\
    \hspace*{1cm}Using continuity of the length functional for any $\varepsilon > 0$ there exists $\delta > 0$ such that for any closed curve $\gamma$ with $d(\gamma_{1}, \gamma) < \delta$ it implies that $L(\gamma) < L_{1} + \varepsilon$. We can pick a simple closed curve $\gamma_{z}^{\varepsilon}$ such that $d(\gamma_{1}, \gamma_{z}^{\varepsilon}) < \delta$. Denote the disc bounded by $\gamma_{z}^{\varepsilon}$ by $C$. Using Theorem 1.8 in~\cite{LNRDisc} we obtain homotopy between any two subarcs of $C$ intersecting only at their endpoints with the lengths bounded by
    \begin{equation*}
        |\partial C| + 686\sqrt{\Area(C)} + 2\diam(C) \leq L_{1} + \varepsilon + 2744\sqrt{2}\frac{A}{L_{1}} + 4d_{z} + L_{1} \leq (2804\sqrt{2} + 64)\frac{A}{L_{1}} + \varepsilon.
    \end{equation*}
    Taking $\varepsilon \rightarrow 0$ we obtain a homotopy between $\gamma_{x}$ and $\gamma_{y}$ inside $M_{z}$ with lengths bounded by $(2804\sqrt{2} + 64)\frac{A}{L_{1}}$. So, the map $f$ is formed by concatenation of these three homotopies and homotopy between the point curves $x$ and $y$ (see Fig.~\ref{short_starfish1}).\\
    \hspace*{1cm}Let's consider the sweep-out of the round sphere by the intersection of it with planes parallel to the $xy$-plane. Then $f$ corresponds to image of a degree $1$ map of this sweep-out. Thus, $f$ is non-contractible. Thus we obtain a simple closed geodesic $\gamma_{2}$ distinct from $\gamma_{1}$ with length $L_{2}$ and
    \begin{equation*}
        L_{1} L_{2} \leq (2804\sqrt{2} + 64)\frac{A}{L_{1}} L_{1} = (2804\sqrt{2} + 64) A.
    \end{equation*}
    The theorem is proved. \hfill\qedsymbol{}
\begin{theorem}
    Suppose $d_{x}$ or $d_{y} < \frac{16\sqrt{2} A}{L_{1}}$. Then there exists a simple closed geodesic $\gamma_{2}$ distinct from $\gamma_{1}$ with length $L_{2}$ such that 
    \begin{equation*}
        L_{1} L_{2} \leq 64 A.
    \end{equation*}
\end{theorem}
\textbf{Proof.}
    \begin{figure}
        \centering
        \includegraphics[scale=0.5]{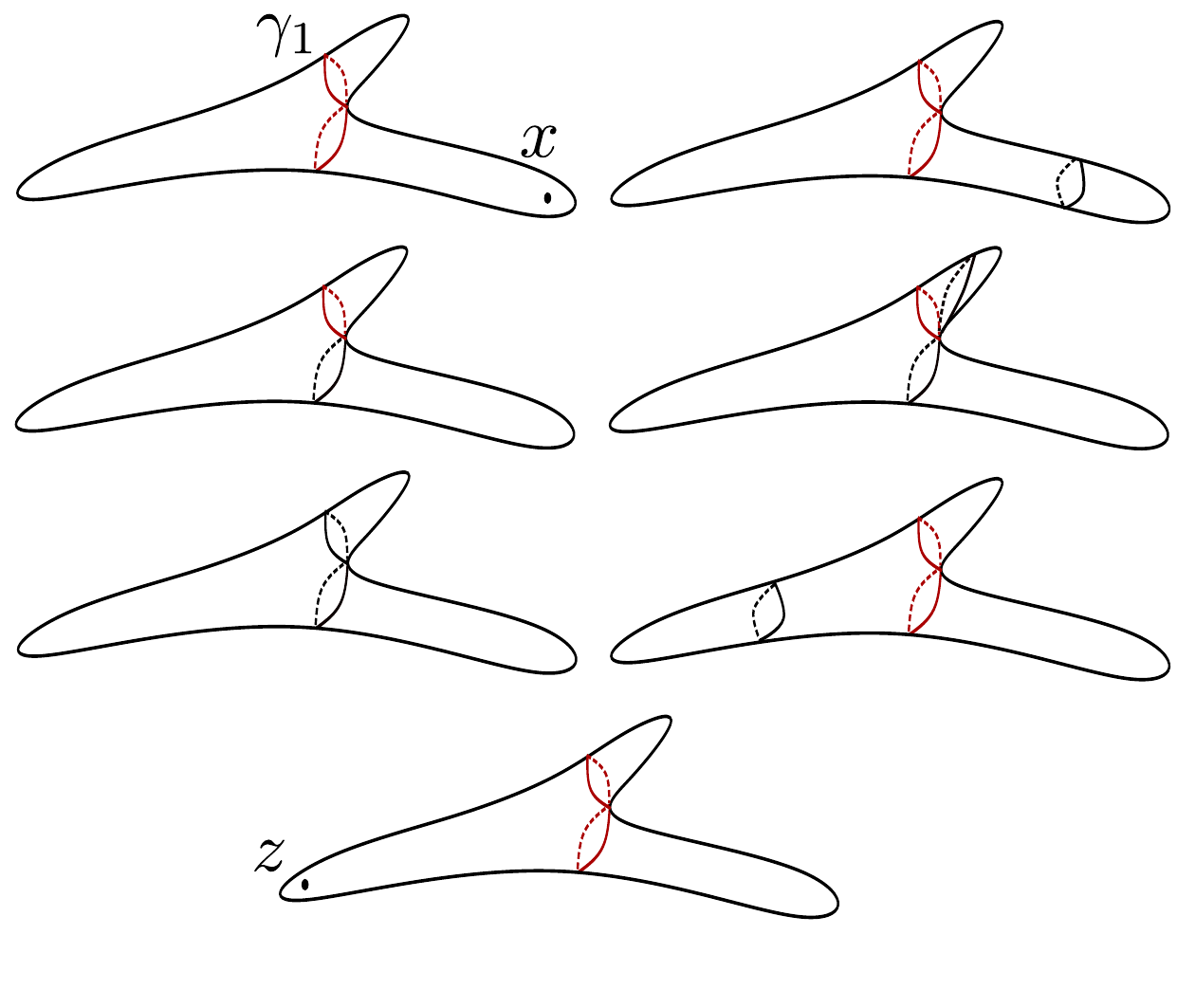}
        \caption{Sweep-out $\phi$}
        \label{short_starfish2}
    \end{figure}
    Without loss of generality we suppose $d_{y} < \frac{16\sqrt{2} A}{L_{1}}$. As in the proof of Theorem 2.4 we will construct a non-contractible map $\phi \: S^1 \rightarrow \Lambda_{L} S^2$. Repeating the argument from the previous proof we suppose we have length-decreasing homotopies of $\gamma_{x}$ with the point curve $x$ inside $M_{x}$ and $\gamma_{y}$ with the point curve $y$ inside $M_{y}$. Now we want to construct a homotopy between $\gamma_{x}$ and $\gamma_{1}$. Let's denote the point of self-intersection of $\gamma_{1}$ by $v$. Connect points $y$ and $v$ by a minimizing geodesic $\tau$. It will lie inside $M_{y}$ and $L(\tau) < \frac{16\sqrt{2} A}{L_{1}}$. We homotope $\gamma_{x}$ to the curve $\gamma_{x}^{} \cup \tau \cup -\tau$. Then we homotope $\gamma_{x} \cup \tau \cup -\tau$ to $\gamma_{1}$ using the homotopy of $\gamma_{y}$ to $y$ (see Fig.~\ref{short_starfish2}). Notice that the lengths of the curves in the obtained homotopy are not greater than
    \begin{equation*}
        L_{1} + \frac{32\sqrt{2} A}{L_{1}} \leq \frac{64\sqrt{2} A}{L_{1}}.
    \end{equation*}
    Again, repeating the same argument with the Birkhoff curve shortening process, we obtain a homotopy of $\gamma_{1}$ with the point curve $z$ inside $M_{z}$ with the lengths not greater than $L_{1}$. The map $\phi$ is non-contractible by the same argument as in the previous proof. So, by the Birkhoff min-max method we obtain a simple closed geodesic $\gamma_{2}$ distinct from $\gamma_{1}$ of length $L_{2}$ not greater than
    \begin{equation*}
        \frac{64\sqrt{2} A}{L_{1}}.
    \end{equation*}
    Then
    \begin{equation*}
        L_{1} L_{2} \leq 64 \sqrt{2} A.
    \end{equation*}
    The theorem is proved. \hfill\qedsymbol{}

\vspace{0.5cm} 
\noindent Department of Mathematics, University of Toronto, Toronto, Canada\\
\textit{E-mail address}: \texttt{talant.talipov@mail.utoronto.ca}


\begin{thebibliography}{99}

\bibitem{LectureNotes} Y. Liokumovich, "Geometric inequalities", [Lecture notes], University of Toronto, 2020. [Online]. \url{https://www.math.toronto.edu/ylio/Geometric_inequalities/}

\bibitem{Gromov} M. Gromov, "Filling Riemannian manifolds", \textit{J. of Differential Geom.}, vol. 18, 1983, 1--147.

\bibitem{Pu} P. M. Pu, "Some inequalities in certain nonorientable Riemannian manifolds", \textit{Pacific J. Math.}, vol. 2, 1952, 55--71.

\bibitem{CrokeKatz} C. B. Croke and M. Katz, "Universal volume bounds in Riemannian manifolds", \textit{Surv. in Differ. Geom.}, vol. 8, 2002, 109--137.

\bibitem{Croke} C. B. Croke, "Area and the length of the shortest closed geodesic", \textit{J. of Differential Geom.}, vol. 27, 1988, 1--21.

\bibitem{Rotman} R. Rotman, "The length of a shortest closed geodesic and the area of a $2$-dimensional sphere", \textit{Proc. Amer. Math. Soc.}, vol. 134, 2006, 3041--3047.

\bibitem{LioNabRot} Y. Liokumovich, A. Nabutovsky, and R. Rotman, "Lengths of three simple periodic geodesics on a Riemannian $2$-sphere", \textit{Math. Ann.}, vol. 367, no. 1--2, 2017, 831--855.

\bibitem{MarquesNevesWillmore} F. C. Marques, A. Neves, "Min-max theory and the Willmore conjecture", \textit{Ann. of Math.}, vol. 179, 2014, no. 2, 683--782.

\bibitem{ChodoshMantoulidis} O. Chodosh, C. Mantoulidis, "The p-widths of a surface," \textit{Publ. Math. IHES}, vol. 137, 2023, 245--342.

\bibitem{Liokumovich} Y. Liokumovich, "Families of short cycles on Riemannian surfaces," \textit{Duke Math. J.}, vol. 165, 2016, no. 7, 1363--1379.

\bibitem{MarquesNeves} F. C. Marques, A. Neves, "Existence of infinitely many minimal hypersurfaces in positive Ricci curvature", \textit{Invent. Math.}, vol. 209, 2017, 577--616.

\bibitem{ColdingMinicozzi} T. H. Colding, W. P. Minicozzi II, "A Course in Minimal Surfaces", \textit{Graduate Studies in Mathematics}, vol. 121, Amer. Math. Soc., Providence, RI, 2011.

\bibitem{Bangert} V. Bangert, "Closed geodesics on complete surfaces", \textit{Math. Ann.}, vol. 251, 1980, 83--96.

\bibitem{ChambLio} G. R. Chambers, Y. Liokumovich, "Optimal sweepouts of a Riemannian $2$-sphere", \textit{J. Eur. Math. Soc.}, vol. 21, 2019, 1361--1377.

\bibitem{LNRDisc} Y. Liokumovich, A. Nabutovsky, and R. Rotman, "Contracting the boundary of a Riemannian $2$-disc", \textit{Geom. Funct. Anal.}, vol. 25, 2015, 1543--1574.

\end{thebibliography}
\end{document}